\documentclass{l4dc2025}

\usepackage{hyperref}
\usepackage{nicefrac}
\usepackage{tikz}
\usetikzlibrary[pgfplots.groupplots]
\usepackage{pgfplots}
\pgfplotsset{compat=1.18}
\usepackage{mathtools}
\usepackage{times}
\usepackage{comment}
\usepackage{caption} %
\newcommand{\diff}[1][]{\mathrm{d}#1}
\newcommand{\dt}{\diff{t}}

\newcommand{\td}[2]{\frac{\diff #1}{\diff #2}}
\newcommand{\argmin}{\displaystyle\mathop{\mathrm{argmin}}}
\newcommand{\minimize}{\displaystyle\mathop{\mathrm{minimize}}}

\newcommand{\vph}{\boldsymbol{\varphi}}

\newcommand{\vps}{\boldsymbol{\psi}}

\newcommand{\vb}{\boldsymbol b}

\newcommand{\ve}{\boldsymbol e}
\newcommand{\vf}{\boldsymbol f}
\newcommand{\vg}{\boldsymbol g}

\newcommand{\vq}{\boldsymbol q}

\newcommand{\vu}{\boldsymbol u}

\newcommand{\vx}{\boldsymbol x}

\newcommand{\vz}{\boldsymbol z}

\newcommand{\vA}{\boldsymbol A}
\newcommand{\vB}{\boldsymbol B}
\newcommand{\vC}{\boldsymbol C}

\newcommand{\vG}{\boldsymbol G}

\newcommand{\vI}{\boldsymbol I}

\newcommand{\vL}{\boldsymbol L}

\newcommand{\nx}{n_{\vx}}
\newcommand{\nU}{n_{\vu}}
\newcommand{\nG}{n_{\vg}}
\newcommand{\nz}{n_{\vz}}

\newcommand{\Plowhat}{\hat{\mathcal{P}}_{\text{lower}}}
\newcommand{\Phighhat}{\hat{\mathcal{P}}_{\text{upper}}}
\newcommand{\vxR}{\vx_T}
\newcommand{\vxL}{\vx_0}
\title[gEDMD Trajectory Optimization]{Koopman Based Trajectory Optimization with Mixed Boundaries}

\coltauthor{
 \Name{Mohamed Abou-Taleb}~\textsuperscript{$\star$} \Email{mohamed.abou-taleb@inm.uni-stuttgart.de}\\
 \Name{Maximilian Raff}~\textsuperscript{$\star$}  \Email{maximilian.raff@inm.uni-stuttgart.de}\\
 \Name{Kathrin Fla{\ss}kamp}~\textsuperscript{$\dagger$}  \Email{kathrin.flasskamp@uni-saarland.de}\\
 \Name{C.~David Remy}~\textsuperscript{$\star$}  \Email{david.remy@inm.uni-stuttgart.de}\\
 \addr \textsuperscript{$\star$}~Institute for Nonlinear Mechanics, University of Stuttgart, Germany\\ \textsuperscript{$\dagger$}~Systems Modeling and Simulation, Saarland University, Germany}

\begin{document}

\maketitle

\begin{abstract}%
Trajectory optimization is a widely used tool in the design and control of dynamical systems.
Typically, not only nonlinear dynamics, but also couplings of the initial and final condition through implicit boundary constraints render the optimization problem non-convex.
This paper investigates how the Koopman operator framework can be utilized to solve trajectory optimization problems in a (partially) convex fashion.
While the Koopman operator has already been successfully employed in model predictive control, the challenge of addressing mixed boundary constraints within the Koopman framework has remained an open question.
We first address this issue by explaining why a complete convexification of the problem is not possible.
Secondly, we propose a method where we transform the trajectory optimization problem into a bilevel problem in which we are then able to convexify the high-dimensional lower-level problem.
This separation yields a low-dimensional upper-level problem, which could be exploited in global
optimization algorithms.
Lastly, we demonstrate the effectiveness of the method on two example systems: the mathematical pendulum and the compass-gait walker.
\end{abstract}


\section{Introduction}
Trajectory optimization problems with mixed boundary conditions (MBCs)  that couple initial and terminal conditions through implicit constraints, are of significant interest across various applications. 
For example, in the field of legged locomotion, trajectory optimization is used as a tool for the design and control of legged robots \citep{Wensing2024}.
Here, MBCs emerge because we want to find optimally actuated periodic solutions for which an initial condition can only be implicitly defined by a periodicity constraint. 
Furthermore, the non-smooth nature of the contact dynamics models render the boundary constraints nonlinear, which adds to the computational challenge.

Due to the nonlinear dynamics, the nonlinear MBCs, and the unknown period time, trajectory optimization problems are typically non-convex. 
Finding a (globally) optimal solution is thus difficult if not impossible, as numerical solvers are of local nature.
Therefore, the quality of the obtained (locally) optimal solution heavily depends on the initial guess provided by the user. 

Figure~\ref{fig:Harmonic Oscillator Optimization} illustrates a simple trajectory optimization problem with MBCs.
Although the cost function is convex in the control input~$u(\cdot)$, the dynamics are linear, and the MBCs are affine, incorporating the period time~$T$ as a decision variable makes the optimization problem non-convex. However, for any fixed $T$, the optimization problem remains convex, allowing us to uniquely determine the optimal solution for that fixed period. 
Exploiting this structure, we can formulate a secondary problem that optimizes over~$T$ in a one-dimensional space. 
In this work, we focus on approximating a broad class of trajectory optimization problems with MBCs, paricularly in the context of nonlinear dynamics, by leveraging Koopman operator theory.
This approach enables us to achieve a structure similar to the introductory example presented in Figure~\ref{fig:Harmonic Oscillator Optimization}, where the subproblem involves linear dynamics, making it convex.    
\begin{figure}
    \begin{minipage}{0.65\textwidth}
    \raggedright
        \begin{flalign*} 
        &\min_{T}
        \begin{cases}
            \begin{array}{cl}
                    \displaystyle\min_{\vx(\cdot), \vu(\cdot)} & c := \displaystyle \int\limits_{0}^{T} u(t)^{2}\dt \\[5mm]
                    \text{s.t.} & \dot{\vx}(t) = \vA\vx(t) + \vB u(t), ~~ \forall t \in [0, T],\\[1mm]
                    &  \underbrace{\begin{bmatrix}
                        \vx(T) - \vx(0)\\
                        \vx(0) - \begin{bmatrix}
                            a \\ 0
                        \end{bmatrix}
                    \end{bmatrix}=\vec{0}}_{\text{mixed boundary conditions}}, \quad\begin{array}{l}
                        \text{(periodicity)}\\
                        \text{(amplitude)}\\
                        \text{(anchor)}
                    \end{array}
            \end{array}
        \end{cases}
        \end{flalign*}
        \caption*{(a)}
 	\label{Harmonic oscillator trajectory optimization problem}
    \end{minipage}\hspace{-3mm}
    \begin{minipage}{0.3\textwidth}
        \begin{tikzpicture}

\begin{axis}[%
width=1\textwidth,
height=0.9\textwidth,
scale only axis,
xmin=0.7,
xmax=5,
xlabel={$T \left[2\pi\right]$},
ymin=-0.01,
ymax=0.24,
xmajorgrids,
ymajorgrids,
legend style={legend cell align=left, align=left}
]
\addplot [color=black, style = thick]
  table[row sep=crcr]{%
0.01	0.00861285444676231\\
0.04125	0.0355278020231794\\
0.0725	0.0624394259224316\\
0.10375	0.0893355124625623\\
0.135	0.11618713137021\\
0.16625	0.142939843130163\\
0.1975	0.169503224800883\\
0.22875	0.19573827107352\\
0.26	0.221442387666825\\
0.29125	0.246332047418995\\
0.3225	0.270023858956182\\
0.35375	0.292015970228349\\
0.385	0.311673558056585\\
0.41625	0.328224680005853\\
0.4475	0.340775670199616\\
0.47875	0.348357512888028\\
0.51	0.350014145352849\\
0.54125	0.34493743870259\\
0.5725	0.332639166031114\\
0.60375	0.313128596876217\\
0.635	0.287043382863331\\
0.66625	0.25567566545968\\
0.6975	0.220857745769885\\
0.72875	0.184720019114242\\
0.76	0.149385236822967\\
0.79125	0.116687171904336\\
0.8225	0.0879847023922212\\
0.85375	0.0640978275093937\\
0.885	0.0453479425787155\\
0.91625	0.0316608202322313\\
0.9475	0.0226897026180561\\
0.97875	0.0179284219548441\\
1.01	0.0167997737931162\\
1.04125	0.0187161074168595\\
1.0725	0.0231156579249754\\
1.10375	0.0294805548355749\\
1.135	0.0373424116500478\\
1.16625	0.0462803221402281\\
1.1975	0.0559148004560902\\
1.22875	0.0659000935038089\\
1.26	0.0759164842781776\\
1.29125	0.0856636799855443\\
1.3225	0.0948560732934631\\
1.35375	0.10322049950597\\
1.385	0.11049700669153\\
1.41625	0.116443031446047\\
1.4475	0.120841155939304\\
1.47875	0.123510251370933\\
1.51	0.124319261000952\\
1.54125	0.123202175609501\\
1.5725	0.120172026409191\\
1.60375	0.115331178387464\\
1.635	0.108875116274699\\
1.66625	0.101087497411248\\
1.6975	0.0923255510252585\\
1.72875	0.0829967153112676\\
1.76	0.0735292522969904\\
1.79125	0.0643408967630022\\
1.8225	0.0558099498357148\\
1.85375	0.0482525254761915\\
1.885	0.0419081750898169\\
1.91625	0.0369343582060405\\
1.9475	0.0334087045535015\\
1.97875	0.0313370506768298\\
2.01	0.0306649101931638\\
2.04125	0.0312902308981856\\
2.0725	0.0330757920089286\\
2.10375	0.0358601934548461\\
2.135	0.0394669388371924\\
2.16625	0.0437115392267708\\
2.1975	0.0484068488392763\\
2.22875	0.0533670025399219\\
2.26	0.0584103891590189\\
2.29125	0.0633620934103909\\
2.3225	0.0680561952712411\\
2.35375	0.0723382426676561\\
2.385	0.076068117489006\\
2.41625	0.079123398323896\\
2.4475	0.0814031874361183\\
2.47875	0.082832219549842\\
2.51	0.083364918117451\\
2.54125	0.0829889317285111\\
2.5725	0.0817275973901923\\
2.60375	0.0796407689791685\\
2.635	0.076823541820354\\
2.66625	0.07340260421926\\
2.6975	0.0695302341603344\\
2.72875	0.0653762873172791\\
2.76	0.0611188252687044\\
2.79125	0.0569342430985334\\
2.8225	0.0529878250685254\\
2.85375	0.049425571895218\\
2.885	0.0463679278744821\\
2.91625	0.0439057445074839\\
2.9475	0.0420985146049123\\
2.97875	0.0409746545609175\\
3.01	0.0405334390098732\\
3.04125	0.0407481125037448\\
3.0725	0.0415697068759581\\
3.10375	0.0429311573744921\\
3.135	0.0447514083341502\\
3.16625	0.0469393058984487\\
3.1975	0.0493971733272962\\
3.22875	0.0520240431365322\\
3.26	0.0547185751550795\\
3.29125	0.0573817202284676\\
3.3225	0.059919197883928\\
3.35375	0.0622438461379191\\
3.385	0.064277876652319\\
3.41625	0.0659550328959793\\
3.4475	0.0672226076558541\\
3.47875	0.0680432347467291\\
3.51	0.0683963344569932\\
3.54125	0.0682790699084834\\
3.5725	0.0677066683928064\\
3.60375	0.0667119823619088\\
3.635	0.0653442103035889\\
3.66625	0.0636667650021853\\
3.6975	0.0617543576758625\\
3.72875	0.0596894493365743\\
3.76	0.0575582918260264\\
3.79125	0.0554468279190136\\
3.8225	0.0534367343381345\\
3.85375	0.0516018712959445\\
3.885	0.050005351620264\\
3.91625	0.0486973712790119\\
3.9475	0.0477138636505998\\
3.97875	0.0470759645216032\\
4.01	0.0467902132860917\\
4.04125	0.046849373881354\\
4.0725	0.0472337380441027\\
4.10375	0.0479127713416304\\
4.135	0.0488469746409957\\
4.16625	0.0499898547753982\\
4.1975	0.0512899228659891\\
4.22875	0.0526926626810648\\
4.26	0.0541424314415399\\
4.29125	0.0555842697909131\\
4.3225	0.0569656055892214\\
4.35375	0.0582378380643899\\
4.385	0.059357785708009\\
4.41625	0.0602889747280643\\
4.4475	0.0610027368419094\\
4.47875	0.0614790778845981\\
4.51	0.0617072742899099\\
4.54125	0.0616861548854187\\
4.5725	0.0614240319870663\\
4.60375	0.0609382590107408\\
4.635	0.0602544112111363\\
4.66625	0.059405110073761\\
4.6975	0.0584285377002947\\
4.72875	0.0573667119975445\\
4.76	0.0562636132674254\\
4.79125	0.055163265098737\\
4.8225	0.0541078754015946\\
4.85375	0.0531361368929513\\
4.885	0.0522817708944195\\
4.91625	0.0515723767707359\\
4.9475	0.0510286241965339\\
4.97875	0.0506638001344665\\
};
\addlegendentry{$\text{c}^\text{*}(T)$}

\addplot[only marks, mark=*, mark options={}, mark size=2pt, color=orange, fill=orange] table[row sep=crcr]{%
x	y\\
2.01007338586448	0.0306649065565302\\
};
\addlegendentry{$c^\ast(T_{\text{local}})$ local optim.}

\addplot[only marks, mark=*, mark options={}, mark size=2pt, color=red, fill=red] table[row sep=crcr]{%
x	y\\
1.00503788267719	0.0167606364369878\\
};
\addlegendentry{$c^\ast(T^\ast)$ global optim.}

\end{axis}

\end{tikzpicture}%
        \caption*{\quad\quad(b)}
 	\label{fig:Multi Search Harmonic Oscillator}
    \end{minipage}
    \caption{Trajectory optimization problem of the harmonic oscillator shown in (a) and its solution family for fixed periods $T$ in (b). The oscillator ($\vA=\left[\begin{smallmatrix}
        1&0\\-1&-0.2
    \end{smallmatrix}\right]$, $\vB=\left[\begin{smallmatrix}
        0\\1
    \end{smallmatrix}\right]$) operates at a desired amplitude of $a = 30$ while maintaining periodicity, satisfying mixed boundary conditions. These conditions, being affine in boundary states with linear dynamics, ensure a unique solution for fixed~$T$. However, allowing $T$ to vary introduces non-convexity, leading to multiple local minima as illustrated in (b).}
    \label{fig:Harmonic Oscillator Optimization}
\end{figure}
In optimal control, convexity is a fundamental property for the efficient synthesis of control policies and system dynamics~\citep{Boyd2004}.
Consequently, the convex formulation of optimization problems has attracted significant attention~\citep{Horst1996,Scherer2006}.
A prominent method to achieving this is (exact) convexification, which transforms the original optimization problem into a convex one. 
A well-known technique in this context is semidefinite relaxation, which expands non-convex constraint spaces into convex ones, often by introducing slack variables to "lift" the problem into a higher-dimensional space~\citep{Acikmese2011, Lasserre2001}.
If the relaxation is proven to be tight, equivalence with the original problem is guaranteed.
Another approach to (exact) convexification, explored in this paper, leverages Koopman operator theory~\citep{Koopman1931}.
This method lifts the nonlinear dynamics into a higher-dimensional space, where they are represented by linear dynamics, thereby enabling a convex problem formulation. 
However, since the Koopman operator is typically infinite-dimensional and challenging to compute, approximate methods, such as the extended dynamic mode decomposition (EDMD), are employed to lift the dynamics into finite-dimensional spaces \citep{Williams2015, Mauroy2020,Williams2016,Iacob2024,Proctor2018}. 
This approach has been successfully implemented in model predictive control (MPC), where it has enabled the optimization problems to be solved efficiently and in real time \citep{Bruder2019a,Korda2018,Kanai2022,Schaller2023}. 
While those works apply the Koopman framework to MPC, they do not consider mixed boundary constraints, since in MPC the prediction horizon and the initial condition are fixed.
In their review of the Koopman operator in robot learning, \cite{Shi2024} discuss the challenge of incorporating constraints within the Koopman framework, which remains an open problem. 

In this work, we investigate how the Koopman framework may be leveraged to simplify trajectory optimization problems with mixed boundary constraints. 
We explain why a complete convexification is not possible. Furthermore, we propose a method where the original problem is approximated by formulating a bilevel optimization problem, encompassing both an upper- and a lower-level problem. This structure allows us to solve the lower-level with fixed boundaries and fixed terminal time as a convex optimization problem which is similar to the MPC subproblem. We present three different ways to formulate the boundary constraints in the lower-level. The upper-level problem then optimizes the boundary values and the terminal time. 
That is, while we still need to solve a nonlinear program, it has a significantly reduced dimensionality.
The efficacy of the presented method is investigated along two examples from periodic trajectory optimization.

\section{Theory}\label{sec: Theory}
\subsection{Problem Definition}
We aim to solve a trajectory optimization problem with mixed boundary constraints (MBCs) of the following form:\\[-6mm]
\begin{figure}[h!]
\centering

\begin{minipage}[t]{0.25\textwidth}
\begin{flushright}
\begin{equation*}
    \mathcal{P}:\left\{\vphantom{\begin{array}{c} \\[18mm] \end{array}}\right.
\end{equation*}
\end{flushright}
\end{minipage}\hspace{-0.25\textwidth}%
\begin{minipage}[t]{1\textwidth}
\begin{subequations}\label{eqn:P}
\begin{flushleft}
\begin{alignat}{2}
    &\minimize_{\vx(\cdot), \vu(\cdot), T} &\quad& \displaystyle c\big( \vx(\cdot), \vu(\cdot), T\big) \notag\\
    &\text{subject to} &\quad& \dot{\vx}(t) = \vf\big(\vx(t)\big) + \vG\big(\vx(t)\big)\vu(t), \quad \forall t \in [0, T], \label{eq:dynamics}\\[2mm]
    & &\quad& \vb\big(\vx(0), \vx(T), T\big) = \vec{0}.\label{eq:MBC}
\end{alignat}
\end{flushleft}
\end{subequations}
\end{minipage}

\end{figure}\\
The trajectories $\vx(t) \in \mathbb{R}^{\nx}$ and $\vu(t) \in \mathbb{R}^{\nU}$ provide the state and input of a control-affine dynamical system with the dynamics $\vf: \mathbb{R}^{\nx} \to \mathbb{R}^{\nx}$ and $\vG: \mathbb{R}^{\nx} \to \mathbb{R}^{\nx\times\nU}$, which is evaluated from time $t=0$ to the terminal time $T$.
Cost is given with the Meyer term $c : \mathbb{R}^{\nx} \times \mathbb{R}^{\nU} \times \mathbb{R} \to \mathbb{R}$, which is assumed to be jointly convex in the state and input. 
The MBCs are defined in an implicit form with the function $\vb : \mathbb{R}^{\nx} \times \mathbb{R}^{\nx} \times \mathbb{R}\to \mathbb{R}^{\nG}$. 
We assume that all functions introduced above are smooth.

Solving $\mathcal{P}$ is inherently challenging due to its pronounced non-convexity, stemming from the variable terminal time $T$, the nonlinear dynamics~\eqref{eq:dynamics} and the nonlinear MBCs~\eqref{eq:MBC}.

\subsection{Koopman Generator Surrogate Modeling}
To deal with the nonlinearity in the dynamics~\eqref{eq:dynamics}, we can lift the dynamical system into an infinite-dimensional space, in which the dynamics are linear~\citep{Koopman1931}.
In the following, we briefly summarize, how we obtain a finite-dimensional approximation of these lifted dynamics.

Let us initially limit ourselves to the time-autonomous case, in which the input is equal to zero ($\vu\equiv \boldsymbol{0}$), such that the dynamics are governed only by the differential equation $\dot{\vx}(t) = \vf(\vx(t))$. 
Let $\vps: \mathbb{R}^{\nx} \to \mathbb{R}^{\nz}$, with $\nz>\nx$, be a nonlinear function  
living itself in a Banach-space, which we denote by~$\mathcal{F}$. In the following, the terms \textit{observable} and \textit{lifting function} will be used synonymously to refer to~$\vps$.
The family of Koopman operators $\mathcal{K}_{t}: \mathcal{F} \to \mathcal{F}$ parameterized by the time $t$ is defined as~$
	(\mathcal{K}_{t}\vps)(\vx_0) = \vps \circ \vph_t(\vx_0) = \vps(\vx(t))$, where $\vph_{t}(\vx_{0})$ denotes the solution of the system~\eqref{eq:dynamics} with the initial condition $\vx_{0}$.
The Koopman operator is linear but infinite dimensional, as it acts on functions. 
Moreover, in the context of continuous-time flows, it is possible to define the infinitesimal Koopman generator as~$\mathcal{L}\vps = \lim_{t\to 0} \nicefrac{1}{t}(\mathcal{K}_{t}\vps-\vps)$;
i.e., the time derivative of the lifting function along solutions, also known as the Lie derivative:
\begin{equation}
	\frac{\mathrm{d}}{\mathrm{dt}}\vps(\vx(t)) = (\mathcal{L}\vps)(\vx(t)) = \nabla_{\vx}\vps(\vx(t))\cdot\vf(\vx(t)).
\end{equation}
As in the case of the Koopman operator, the Koopman generator is linear and infinite-dimensional. 

When we now re-introduce  control inputs $\vu(t)$,  \cite{Peitz2020a} state in Theorem 3.2 that the Koopman generators inherit the property of control-affinity. 
Consequently, the Koopman generator $\mathcal{L}_{\vu(t)}$ for a known input function $\vu(t)$ may be expressed as
$\mathcal{L}_{\vu(t)} = \mathcal{L}_{0} + \sum_{i = 0}^{n_{\vu}}u_i(t)(\mathcal{L}_{i} - \mathcal{L}_{0})$,
where $\vu(t) = \sum_{i=1}^{n_{\vu}}u_{i}(t)\ve_{i}$ with the $ i $th canonical basis vector $\ve_{i}$, $\mathcal{L}_{0}$ is the Koopman generator for $\vu \equiv 0$ and each $\mathcal{L}_{i}$ denotes the Koopman generator for each basis element of $\vu(t)$. Thus, for control affine systems, the Lie derivative of the observables yields
\begin{equation}\label{eqn:bilinearGenerator}
	\td{}{t}\vps(\vx(t)) = (\mathcal{L}_{\vu(t)}\vps)(\vx(t))
	= (\mathcal{L}_{0}\vps)(\vx(t)) + \sum_{i = 0}^{n_{\vu}}u_i(t)(\mathcal{L}_{i} - \mathcal{L}_{0})\vps(\vx(t)),
\end{equation}
which is bilinear in the lifted state $\vps(\vx)$ and the input $\vu(t)$.

Here, we focus on approximating the Koopman generator rather than the Koopman operator, because we aim to ultimately optimize over the terminal time $T$.
When discretizing the optimization problem~$\mathcal{P}$ with a fixed number of discrete time steps, changes in $T$ would directly affect the step size and would thus require a repeated identification of the Koopman operator.
In contrast, the Koopman generator allows us to use a continuous-time dynamics description, which is identified once.
This can then be applied to discretizations over arbitrary grids. 

To obtain a finite-dimensional approximation of the Koopman generator, we use the \textit{generator Extended Dynamic Mode Decomposition} (gEDMD), proposed by \cite{Klus2020}, which approximates the generator as:
\begin{equation}\label{eq:gEDMD}
	\td{}{t}\vps(x(t)) \approx \vL_{0}\vps(\vx(t)) + \sum_{i = 1}^{n_{\vu}}u_{i}(t)(\vL_{i}-\vL_{0})\vps(\vx(t)),
\end{equation}
where $\vL_{i}\in \mathbb{R}^{\nz\times\nz}$, $i \in\{ 0,\dots,\nU\}$, are finite dimensional matrices identified from data.

\subsection{Koopman-Based Trajectory Optimization}
\label{sec:bilevel}
Exact convexification of~$\mathcal{P}$ using Koopman operator theory would be possible if the terminal time~$T$ and the initial state $\vec{x}(0)$ are explicitly defined by the MBCs~\eqref{eq:MBC}.
In this case, we could simply lift the initial state with the observables to yield a lifted initial state  $\vz(0) = \vps(\vx(0))$.
This does not only ensure that the resulting solution corresponds to a trajectory starting at $\vec{x}(0)$, but also that the lifted trajectory is bound to the manifold~$\mathcal{M}$, defined by the observable functions:
\begin{equation}
    \mathcal{M} \coloneqq \{\vz \in \mathbb{R}^{\nz} \, \lvert \, \vz = \vps(\vx), \, \vx \in \mathbb{R}^{\nx} \}.
\end{equation}
Analogously, this process could be done with the final state $\vec{x}(T)$, or for any other time instant~$\bar{t}\in[0,T]$ along the trajectory.

For general MBCs~\eqref{eq:MBC}, however, the boundary states and terminal time are not explicitly defined.
While the implicit function $\vb$ could be convexified through a suitable choice of observables~$\vps$, yielding a corresponding lifted constraint~$\vb_\text{lift}\big(\vz(0), \vz(T), T\big) = \vec{0}$,
this alone is insufficient to guarantee a valid solution:
the lifted constraint can be satisfied by choices of $\vz(0)$ and $\vz(T)$ that do not lie within the manifold~$\mathcal{M}$.
Therefore, an additional constraint is required to ensure that, for some $\bar{t}$, it holds: $\vz(\bar{t}) \in \mathcal{M}$.
Since $\vps$ is nonlinear, this additional constraint may not be convex, meaning the resulting lifted optimization problem could be non-convex, even if $\vb$ or $\vb_\text{lift}$ are convex.

A similar issue arises with the terminal time $T$ (and other parameters).
One approach to address the issue of an unknown terminal time $T$ (as illustrated in the introductory example in Figure~\ref{fig:Harmonic Oscillator Optimization}) is to scale all derivatives by $\nicefrac{1}{T}$ and evaluate the dynamics over a normalized time interval~$\tau \in [0,1]$.
This allows us to treat $T$ as a parameter, which could be included into the dynamics as an additional state with~$T^\prime(\tau)= 0$, resulting in an augmented lifted state~$\vz(\tau)=\vps\big(\vx(\tau),T(\tau)\big)$.
However, since neither $T(0)$ nor $T(1)$ are explicitly defined, the same problem arises as above:
$\vz(\bar{\tau})$ must be constrained to lie on $\mathcal{M}$ for some time~$\bar{\tau}\in[0,1]$, requiring an additional constraint that may be non-convex.

To address these issues, we will approximate~$\mathcal{P}$ within a bilevel optimization problem~$\hat{\mathcal{P}}$. 
This encompasses a non-convex upper-level optimization problem~$\Phighhat$ and a convex lower-level optimization problem~$\Plowhat$.
This formulation enables the use of algorithms which attempt to find global optima, because $\Phighhat$ optimizes over a space which is low-dimensional, while the convexification of the lower-level problem alleviates the burden of solving an expensive and high-dimensional nonlinear program. 
We define the upper-level problem as
\begin{equation*}
\Phighhat:
\begin{cases}
    \begin{array}{cl}
            \minimize_{\vxL, \vxR, T} & \displaystyle c\big(\vx^{\ast}(\cdot), \vu^{\ast}(\cdot), T\big) \\[5mm]
            \text{subject to} & \left( \vz^{\ast}(\cdot), \vu^{\ast}(\cdot) \right) = \Plowhat(\vxL, \vxR, T),\\
            & \vx^{\ast}(\cdot) = \vC \vz^{\ast}(\cdot),\\
            & \vb\big(\vxL, \vxR, T\big) = 0,
    \end{array}
\end{cases}
\end{equation*}
and the lower-level as
\begin{equation*} 
    \Plowhat\big(\vxL, \vxR, T\big):
    \begin{cases}
        \begin{array}{cl}
             \argmin_{\vz(\cdot), \vu(\cdot)} & \displaystyle (1-w_i) \cdot c\big( \vC \vz(\cdot), \vu(\cdot); T\big) + w_i\cdot \hat{c}_\text{lower}\big(\vz(0), \vz(T); \vxL, \vxR\big) \\[5mm]
            \text{subject to} & \dot{\vz}(t) =  \vL_{0}\vz(t) + \vL_{\vu}(\bar{\vz})\vu(t), \quad\forall t \in [0, T],\\
          & \hat{\vb}_i\big(\vz(0), \vz(T); \vxL, \vxR\big) = \vec{0},
        \end{array}
    \end{cases}
\end{equation*}
where the lift~$\vz = \vps(\vx)$ includes a copy of $\vx$ such that we can recover the original state by applying a linear operation $\vx = \vC \vz$, with $\vC=\left[
        \vI_{\nx \times \nx} ~ \mathbf{0}\right]$.

In the upper-level~$\Phighhat$, we introduce the initial state~$\vxL$ and final state~$\vxR$ as decision variables, together with the terminal time $T$.
We enforce the MBCs on $\vxL$, $\vxR$ and $T$ in a nonlinear fashion, while the cost is evaluated by obtaining optimal state- and input trajectories for fixed boundary values and terminal time in the lower-level~$\Plowhat$.
Note that the upper-level~$\Phighhat$ is a low-dimensional non-convex problem. 

For each iteration of $\Phighhat$, the bilinear dynamics~$\eqref{eq:gEDMD}$ are linearized in the lower-level~$\Plowhat$ around a chosen point~$\bar{\vz}=\vps(\bar{\vx})$ with $\vu = \mathbf{0}$.
This results in the following lifted linear dynamics:
\begin{equation} \label{eq:LTI}
    \dot\vz(t) = \vL_{0}\vz(t) + \vL_{\vu}(\bar{\vz})\vu(t),
\end{equation}
where $\vL_{\vu} \in \mathbb{R}^{\nx\times\nU}$ is the input matrix of the resulting LTI dynamics.
Furthermore, we incorporate an additional cost~$\hat{c}_\text{lower}$ with weights $w_i$ and boundary constraints $\hat{\vb}_i$ into~$\Plowhat$.

The implementations of $\hat{\vb}_i$ must include $\vC\vz(0) = \vxL$ and $\vC\vz(T) = \vxR$, ensuring that the lifted state adheres to the original state constraints.
As discussed above, this alone is insufficient, as~$\vz$ must additionally be constrained to~$\mathcal{M}$. 
A naive approach to ensure this, would be to fully lift both boundary conditions, setting $\vz(0) = \vps(\vxL)$ and $\vz(T) = \vps(\vxR)$.
However, this leads to practical issues if the lifted state~$\vz$ drifts away from $\mathcal{M}$ due to the approximative nature of $\vL$. 
Since the lifted linear dynamics~\eqref{eq:LTI} are not necessarily controllable in directions orthogonal to~$\mathcal{M}$, no input $\vu$ would be able to steer the lifted state back onto $\mathcal{M}$. 

Alternatively, this work explores and compares three different approaches ($i\in\{0,T,\text{soft}\}$) for lifting the boundary conditions. 
For the first two, we simply apply the lift to only one of the boundaries, yielding the two implicit formulations:
\begin{subequations}
    \begin{align}\label{eq:b0}
        \hat{\vb}_0\big(\vz(0), \vz(T); \vxL, \vxR\big) &= \begin{bmatrix}
            \vz(0) - \vps(\vxL) \\
            \vC\vz(T) - \vxR
        \end{bmatrix},
        \\
        \hat{\vb}_T\big(\vz(0), \vz(T); \vxL, \vxR\big) &= \begin{bmatrix}\label{eq:bT}
            \vC\vz(0) - \vxL \\
            \vz(T) - \vps(\vxR)
        \end{bmatrix},
    \end{align}
where the weights~$w_0=w_T=0$ are chosen to exclusively  account for the original cost~$c$.
As a third approach, we implement the requirements $\vz(0) = \vps(\vxL)$ and $\vz(T) = \vps(\vxR)$ as soft constraints, by invoking the additive penalty cost~$\hat{c}_\text{lower}$ with weights $w_{\text{soft}}\in(0,1)$:
    \begin{align}\label{eq:softCost}
        \hat{c}_\text{lower}(\vz(0), \vz(T); \vxL, \vxR) &=  \lVert\vz(0) - \vps(\vxL) \rVert^{2} +  \lVert\vz(T) - \vps(\vxR) \rVert^{2},\\
        \hat{\vb}_{\text{soft}}\big(\vz(0), \vz(T); \vxL, \vxR\big) &= \begin{bmatrix}
            \vC\vz(0) - \vxL \\
            \vC\vz(T) - \vxR
        \end{bmatrix}.\label{eq:bsoft}
    \end{align}
\end{subequations}
Note: Equations~\eqref{eq:b0} and \eqref{eq:bT} will enforce the lifted state to lie on~$\mathcal{M}$ only at the beginning or end of the trajectory, respectively.
Since $\mathcal{M}$ is only approximately Koopman-invariant under $\vL$, the lifted state will drift away from $\mathcal{M}$ throughout the remainder of the trajectory.
Equations~\eqref{eq:softCost} and ~\eqref{eq:bsoft} seek to balance this drift across both boundaries, ideally yielding a solution that is closer to $\mathcal{M}$ overall.
Yet, it may be prone to the ill-conditioning associated with soft constraints \citep{Betts2010}.
Finally, we have to select an appropriate point of linearization for each case. 
Given that the only known points on the solution trajectory are the initial and terminal state, it is reasonable to use these for the linearization.
In the case where we lift the initial constraints $\hat{\vb}_{0}$, we should choose $\bar{\vz} = \vps(\vxL)$.
If we lift the terminal constraint $\hat{\vb}_{T}$, we should choose $\bar{\vz} = \vps(\vxR)$.
In the case of soft constraints, both options are equally valid.
This linearization yields $\vL_{\vu}$ from $\vL_1$ and is conducted repeatedly for each call to the lower-level problem. 

The cost functions $c$ and $\hat{c}_\text{lower}$, as well as convex combinations of these functions, are jointly convex in the decision variables~$\vz(\cdot)$ and $\vu(\cdot)$. 
Furthermore, the constraints $\hat{\vb}_i\big(\vz(0), \vz(T); \vxL, \vxR\big)$, with $i\in\{0,T,\text{soft}\}$, are affine and define a convex feasible set. 
Consequently, $\Plowhat$ is a convex optimization problem.

\section{Examples from Periodic Trajectory Optimization}\label{sec: Examples}
\begin{figure}[t]
    \centering
    \includegraphics{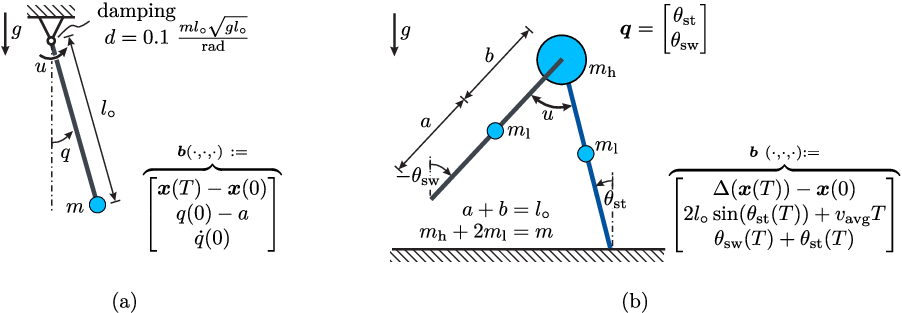}
    \caption{Illustrated are two mechanical systems: (a) a mathematical pendulum and (b) a compass-gait walker. The MBCs $\vb$ include periodicity, operating, and anchor constraints. For (a), the operating constraint is the desired amplitude~$a$, and for (b), it is the desired average forward velocity~$v_\text{avg}$. The system states are $\vx^\top = [\vq^\top~\dot{\vq}^\top]$, with all parameters normalized by gravity~$g$, length~$l_\circ$ and mass~$m$.}
    \label{fig:examples}
\end{figure}

To investigate the efficacy of the presented method, we consider two (periodic) example systems, namely a mathematical pendulum and a compass-gait walker (Figure~\ref{fig:examples}). 
For the numerical evaluation of these problems, all system parameters are normalized with respect to gravity~$g$, length~$l_\circ$ and mass~$m$.
To numerically solve the optimization problems~$\mathcal{P}$ and $\hat{\mathcal{P}}$, the input space is approximated with piecewise constant functions. 
Nonlinear dynamics within $\mathcal{P}$ are approximated with an explicit fourth-order Runge-Kutta integration scheme whereas linear dynamics within~$\Plowhat$ are discretized exactly using matrix exponentials.
In both examples, we utilize the cost function~$c(\vx(\cdot), \vu(\cdot), T) = \int_{0}^{T} u(t)^{2}\dt$.
An optimal solution for $\hat{\mathcal{P}}$ is obtained utilizing MATLAB's \texttt{fmincon} for $\Phighhat$ and \texttt{quadprog} for $\Plowhat$.
To evaluate the solution of our proposed bilevel optimization problem, we also solve the resulting nonlinear program of the original problem $\mathcal{P}$, again using \texttt{fmincon}. Here, we use the solution of $\hat{\mathcal{P}}$ as the initial guess for  $\mathcal{P}$.
For further details on the implementation, including the specific choice of lifting functions, we refer to the actual code that is available on \href{https://github.com/MohamedAbou-Taleb/KoopmanBasedTrajectoryOptimizationWithMixedBoundaries}{GitHub}\footnote{The code can be found at the following GitHub repository: \url{https://github.com/MohamedAbou-Taleb/KoopmanBasedTrajectoryOptimizationWithMixedBoundaries}}. 

\subsection{Mathematical Pendulum}
\begin{figure}[t]
\centering
    \begin{minipage}{0.4\textwidth}
        \centering
        \input{plots/resultsSP/simplePendulumError}
        \caption*{(a)}
    \end{minipage}\hfill
    \begin{minipage}{0.55\textwidth}
        \centering
        \input{plots/resultsSP/simplePendulumPhaseSpace}
        \caption*{(b)}
    \end{minipage}
    \caption{Results of the mathematical pendulum, compared across different choices for the boundary constraints $\hat{\vb}_i$ and weights $w_{\text{soft}}$. The discretization utilizes $101$ points.
    (a) evaluates the similarity of optimal period time $T^\ast$, optimal states $\vx^\ast(\cdot)$ and optimal inputs $\vu^\ast(\cdot)$, as a function of amplitude $a$ between the approximated problem and the original NLP.
    Similarity of trajectories is expressed via the Pearson correlation coefficient (PCC). 
    (b) shows a phase portrait for an amplitude of $a=40^\circ$.
    }
    \label{fig:Simple Pendulum results}
\end{figure}

As a nonlinear extension of the harmonic oscillator trajectory optimization problem presented in the introduction (Figure~\ref{fig:Harmonic Oscillator Optimization}), we consider the simple pendulum (Figure~\ref{fig:examples}a).
Apart from the non-linear dynamics, the optimization problem is identical, seeking periodic trajectories with a desired amplitude $a$ that are anchored at a velocity of $\dot{q}(0)=0$.
This problem was rewritten in the bilevel form~$\hat{\mathcal{P}}$, introduced in Section~\ref{sec:bilevel}.
To convexify the dynamics, we lifted them by identifying a bilinear Koopman generator surrogate model. To this end, we utilized a 12-dimensional Koopman basis consisting of trigonometric functions and polynomials of the state. We collected data by sampling $45000$ uniformly distributed points in the state space, subsequently lifted them, and computed their Lie-derivatives. With this data we solved two least-squares problems to obtain the matrices $\vL_{0}, \, \vL_{1}$ for the bilinear approximation of the Koopman generator.
We solved this problem with all three versions of the lifted boundary conditions $\hat{\vb}_i$ and for three values of $w_{\text{soft}}\in\{0.1, 0.5, 0.9\}$, and we compare the results to the solution obtained from numerically solving the original problem~$\mathcal{P}$~(Figure~\ref{fig:Simple Pendulum results}).
In all cases, we linearized the dynamics about $\bar{\vz} = \vps(\vxL)$ (which, given the periodicity, is equivalent to linearizing about $\bar{\vz} = \vps(\vxR)$).

Figure~\ref{fig:Simple Pendulum results}a shows the similarity of optimal period time $T^\ast$, optimal states $\vx^\ast(\cdot)$ and optimal inputs $\vu^\ast(\cdot)$, as a function of amplitude $a$.
The similarities of the latter are expressed via the Pearson correlation coefficient (PCC), for which a full agreement between original and approximated trajectories would yield a value of $1$.
We observe that $\hat{\vb}_0$ and $\hat{\vb}_T$ lead to the best agreement, while the soft constraint are introducing a larger approximation error, in particular if the relative weight of the two cost components is unbalanced.

For an amplitude of $a=40^{\circ}$, the cost of $\mathcal{P}$ yields $c = 1.48\cdot 10^{-2} (m g l_{\circ})^{2}\sqrt{\nicefrac{l_{\circ}}{g}}$, while the hard constraints with the lift on the initial- and the terminal condition yield a cost of $c = 3.34\cdot 10^{-2} (m g l_{\circ})^{2}\sqrt{\nicefrac{l_{\circ}}{g}}$ and $c = 1.53\cdot 10^{-2} (m g l_{\circ})^{2}\sqrt{\nicefrac{l_{\circ}}{g}}$ as well as $\hat{c}_\text{lower} = 23.88\cdot 10^{-2}$ and $\hat{c}_\text{lower} = 102.21\cdot 10^{-2}$ respectively. The costs in the case of the soft constraints evaluate to $c = 0.22\cdot 10^{-2} (m g l_{\circ})^{2}\sqrt{\nicefrac{l_{\circ}}{g}}$ and $\hat{c}_\text{lower} = 3.013\cdot 10^{-2}$ for $w_{\text{soft}}=0.1$;
$c = 1.02\cdot 10^{-2} (m g l_{\circ})^{2}\sqrt{\nicefrac{l_{\circ}}{g}}$ and $\hat{c}_\text{lower} = 1.34\cdot 10^{-2}$ for $w_{\text{soft}}=0.5$;
as well as $c = 4.00\cdot 10^{-2} (m g l_{\circ})^{2}\sqrt{\nicefrac{l_{\circ}}{g}}$ and $\hat{c}_\text{lower} = 0.84\cdot 10^{-2}$ for $w_{\text{soft}}=0.9$. 
As expected, for the soft constraints, we see a clear trade-off between cost $c$ and accuracy of the solution, as measured by $\hat{c}_\text{lower}$.  
As per our initial hypothesis, enforcing the constraints only at one boundary via the hard constraints, leads to much larger constraint violations at the other end of the lifted trajectory (expressed here via $\hat{c}_\text{lower}$).
However, when comparing this to the results presented in Figure~\ref{fig:Simple Pendulum results}, the better performance of the soft constraints in the lifted space, does not translate to a better agreement in terms of un-lifted trajectories. 

\subsection{Compass-Gait Walker}
\begin{figure}[t]
\centering
    \begin{minipage}{0.4\textwidth}
        \centering
        \definecolor{Hellblau}{rgb}{0.00000,0.74510,1.00000}%
\begin{tikzpicture}

\begin{axis}[%
width=1\textwidth,
height=160pt,
xmin=0,
xmax=2.5,
xlabel={$t\, [\sqrt{l_\circ/g}]$},
ymin=0,
ymax=0.0005,
ylabel={$u^\ast\, [m gl_\circ]$},
ytick={0, 0.0002, 0.0004},
yticklabels={0, 2, 4},
scaled ticks=false,
xmajorgrids,
ymajorgrids,
legend style={legend cell align=left, align=left, at={(1.7in,1.5in)}}
]
\addplot [const plot, color=Hellblau, style = thick]
  table[row sep=crcr]{%
0	0.000143850234455836\\
0.045270582737833	0.000166068314161321\\
0.090541165475666	0.000188023982350596\\
0.135811748213499	0.000209612148872854\\
0.181082330951332	0.000230728374538822\\
0.226352913689165	0.000251269593841907\\
0.271623496426998	0.000271134804370486\\
0.316894079164831	0.000290225722703634\\
0.362164661902664	0.000308447406489739\\
0.407435244640497	0.000325708842333783\\
0.45270582737833	0.000341923499063509\\
0.497976410116163	0.000357009845905224\\
0.543246992853996	0.000370891835076253\\
0.588517575591829	0.000383499348291698\\
0.633788158329662	0.000394768606687552\\
0.679058741067495	0.000404642543679169\\
0.724329323805328	0.000413071140303213\\
0.769599906543161	0.000420011722631695\\
0.814870489280994	0.000425429220897637\\
0.860141072018827	0.000429296390032996\\
0.905411654756661	0.000431593991389814\\
0.950682237494493	0.000432310935494759\\
0.995952820232326	0.000431444385774343\\
1.04122340297016	0.000428999823282872\\
1.08649398570799	0.000424991072566674\\
1.13176456844583	0.000419440288905959\\
1.17703515118366	0.000412377907288829\\
1.22230573392149	0.000403842553590056\\
1.26757631665932	0.000393880918549418\\
1.31284689939716	0.000382547595269989\\
1.35811748213499	0.000369904881084837\\
1.40338806487282	0.000356022544770718\\
1.44865864761066	0.000340977560218358\\
1.49392923034849	0.000324853807800202\\
1.53919981308632	0.000307741744807286\\
1.58447039582416	0.00028973804645617\\
1.62974097856199	0.000270945219093971\\
1.67501156129982	0.000251471187353521\\
1.72028214403765	0.000231428857130858\\
1.76555272677549	0.000210935656372661\\
1.81082330951332	0.000190113055771188\\
1.85609389225115	0.000169086071567996\\
1.90136447498899	0.000147982752764396\\
1.94663505772682	0.000126933655125591\\
1.99190564046465	0.000106071304446071\\
2.03717622320249	8.55296516153866e-05\\
2.08244680594032	6.54435220854914e-05\\
2.12771738867815	4.59480623926533e-05\\
2.17298797141598	2.71781864282213e-05\\
2.21825855415382	9.26802418272796e-06\\
};
\addlegendentry{$\hat{\mathcal{P}}$}
\addplot [const plot,color=red, style = thick]
  table[row sep=crcr]{%
0	0.000141832920892725\\
0.0450366200291482	0.000162768939907747\\
0.0900732400582963	0.00018485485534377\\
0.135109860087444	0.000214465429263984\\
0.180146480116593	0.000230721135560782\\
0.225183100145741	0.000260809471521472\\
0.270219720174889	0.000275221500873814\\
0.315256340204037	0.000297473773687198\\
0.360292960233185	0.000321057510007774\\
0.405329580262333	0.000336052666411826\\
0.450366200291482	0.000349898934984455\\
0.49540282032063	0.000367770010422951\\
0.540439440349778	0.000384477754470639\\
0.585476060378926	0.000395343036499797\\
0.630512680408074	0.00040307293611741\\
0.675549300437222	0.000413099785392266\\
0.720585920466371	0.000426914157278312\\
0.765622540495519	0.000441400711855871\\
0.810659160524667	0.000452905299501409\\
0.855695780553815	0.000460568956721306\\
0.900732400582963	0.00046600453371653\\
0.945769020612111	0.000470609019503599\\
0.99080564064126	0.000473546749874789\\
1.03584226067041	0.000472211999583645\\
1.08087888069956	0.000464626498557597\\
1.1259155007287	0.000451703364738722\\
1.17095212075785	0.000437445377435794\\
1.215988740787	0.00042663153941558\\
1.26102536081615	0.000421377483130687\\
1.3060619808453	0.000419027881577027\\
1.35109860087444	0.000413340562392687\\
1.39613522090359	0.000398808353048822\\
1.44117184093274	0.000375271964701014\\
1.48620846096189	0.000348735231957256\\
1.53124508099104	0.000326450872034293\\
1.57628170102019	0.000309638535907089\\
1.62131832104933	0.000291514534611867\\
1.66635494107848	0.000265305927890607\\
1.71139156110763	0.000235139036600543\\
1.75642818113678	0.000213566877472849\\
1.80146480116593	0.000202185587499682\\
1.84650142119507	0.000183483926339337\\
1.89153804122422	0.000151368281068584\\
1.93657466125337	0.000133054334230465\\
1.98161128128252	0.000126246712580091\\
2.02664790131167	8.27104399673783e-05\\
2.07168452134082	8.04997604443148e-05\\
2.11672114136996	4.0938967219019e-05\\
2.16175776139911	2.92294910697035e-05\\
2.20679438142826	1.20759254647938e-05\\
};
\addlegendentry{$\mathcal{P}$}

\end{axis}
\end{tikzpicture}%
        \caption*{(a)}
    \end{minipage}\hfill
    \begin{minipage}{0.55\textwidth}
        \centering
        \definecolor{Hellblau}{rgb}{0.00000,0.74510,1.00000}%
\begin{tikzpicture}

\begin{axis}[%
width=1\textwidth,
height=160pt,
xmin=-4,
xmax=5,
xlabel style={font=\color{white!15!black}},
xlabel={$\theta_{i}^\ast~\big[{}^\circ\big]$},
ymin=-8,
ymax=10,
ylabel={$\dot{\theta}^\ast_{\text{i}}\big[{}^\circ/\sqrt{l_\circ/g}\big]$},
xmajorgrids,
ymajorgrids,
legend style={legend cell align=left, align=left , at={(1.95in,1.2in)}, legend columns=2}
]
\addplot [color=Hellblau, style = thick]
  table[row sep=crcr]{%
-3.2439995696396	-4.83280576771475\\
-3.44301848041854	-3.96064536640336\\
-3.60269385778515	-3.09519956839813\\
-3.72340407138616	-2.23975953697649\\
-3.80567706242416	-1.39763877127413\\
-3.85019091589511	-0.572153727498115\\
-3.85777355706273	0.233395476979755\\
-3.8294015780606	1.0157441821857\\
-3.7661982047824	1.77168203267624\\
-3.66943041846235	2.49807146003149\\
-3.54050525053798	3.19186570887789\\
-3.38096527350447	3.85012631647917\\
-3.19248331449046	4.47003995829129\\
-2.97685642218921	5.04893457470775\\
-2.73599912154367	5.58429469751374\\
-2.47193599419	6.07377589830187\\
-2.18679362609099	6.51521828526091\\
-1.8827919670208	6.90665897930754\\
-1.5622351495767	7.24634350546578\\
-1.22750181817724	7.53273604068127\\
-0.881035021042831	7.76452846485909\\
-0.525331720431368	7.94064816780202\\
-0.162931978405922	8.06026457086859\\
0.203592122867744	8.12279432853201\\
0.571647765867149	8.1279051815665\\
0.938632825947012	8.07551844027965\\
1.3019473806822	7.96581008301117\\
1.65900523742711	7.79921046199168\\
2.00724541599724	7.57640261556195\\
2.34414352348985	7.29831919265566\\
2.66722295868339	6.96613800230914\\
2.97406588418702	6.58127620774324\\
3.26232390554781	6.14538319122885\\
3.52972839785806	5.66033212246313\\
3.7741004220307	5.12821026951502\\
3.99336017481827	4.55130809751242\\
4.18553591882988	3.93210720610906\\
4.34877234123847	3.27326716235691\\
4.48133829255509	2.57761129089233\\
4.58163385976299	1.84811148829547\\
4.64819673123629	1.08787213307891\\
4.67970781420021	0.300113166982078\\
4.67499606900352	-0.511847572927834\\
4.63304252815205	-1.34461268843565\\
4.55298347187446	-2.19472413474232\\
4.43411273593982	-3.05868262697389\\
4.27588313149868	-3.93296744810427\\
4.07790696085709	-4.81405656612887\\
3.8399556172928	-5.69844696716082\\
3.56195826126595	-6.5826751104265\\
3.2439995696396	-7.46333741091518\\
};
\addlegendentry{sw ($\hat{P}$)}

\addplot [color=Hellblau, style = thick, dashed]
  table[row sep=crcr]{%
3.2439995696396	-5.18834003080495\\
3.0156229176292	-4.90307662935462\\
2.79989049901271	-4.6297132142638\\
2.59626601838627	-4.36814775001355\\
2.40421712663074	-4.11830740857192\\
2.22321416905115	-3.88014547149071\\
2.05272907223172	-3.6536383008014\\
1.8922343663131	-3.43878238663688\\
1.74120233903644	-3.23559147963459\\
1.59910431753247	-3.04409381624243\\
1.46541007346579	-2.8643294450506\\
1.3395873467781	-2.6963476622116\\
1.22110148291156	-2.54020456388744\\
1.10941517803802	-2.39596072348044\\
1.00398832647392	-2.26367900116161\\
0.904277964127183	-2.1434224929138\\
0.809738301503287	-2.03525262595546\\
0.719820839495818	-1.93922740701064\\
0.633974560904074	-1.85539982944353\\
0.551646190358763	-1.78381644478634\\
0.472280515098338	-1.7245161036614\\
0.395320758824752	-1.67752887053643\\
0.320209000679705	-1.64287511616073\\
0.24638663122217	-1.62056479091435\\
0.173294837155936	-1.61059688166698\\
0.100375106453044	-1.61295905409329\\
0.0270697454457705	-1.62762748173238\\
-0.0471776005834929	-1.65456686241463\\
-0.122921431798259	-1.69373062201646\\
-0.200713891656987	-1.74506130484529\\
-0.281104286635957	-1.8084911493096\\
-0.36463863617152	-1.8839428468967\\
-0.45185926106377	-1.97133048186718\\
-0.543304418452979	-2.07056064848628\\
-0.639507991323317	-2.18153374205036\\
-0.740999240306906	-2.30414541943661\\
-0.848302625356469	-2.43828822440912\\
-0.961937704628579	-2.58385337245728\\
-1.08241911767342	-2.7407326895266\\
-1.21025665976307	-2.90882069862962\\
-1.34595545391051	-3.0880168479983\\
-1.49001622683822	-3.27822787416028\\
-1.64293569485018	-3.47937029309267\\
-1.80520706524688	-3.69137301242823\\
-1.97732065860203	-3.91418005756271\\
-2.15976465689387	-4.14775340443844\\
-2.35302598215624	-4.39207591175859\\
-2.5575913099874	-4.64715434542039\\
-2.77394822192934	-4.9130224880425\\
-3.00258650041142	-5.18974432660298\\
-3.2439995696396	-5.47741731139941\\
};
\addlegendentry{st ($\hat{P}$)}
\addplot [color=red, style = thick]
  table[row sep=crcr]{%
-3.22721647058485	-4.83656490068622\\
-3.42559112761699	-3.97378515324227\\
-3.58522668787217	-3.11676512224332\\
-3.70645423518064	-2.26869379543787\\
-3.78974189256126	-1.4324653156816\\
-3.83571203939575	-0.612032979966852\\
-3.84513421508342	0.19002417958067\\
-3.81892704617366	0.969667241234419\\
-3.75816356119819	1.72408146793338\\
-3.66404966398229	2.45018349764893\\
-3.5379397263799	3.14446651684042\\
-3.38133556869607	3.8038926109959\\
-3.1958674823206	4.4257974679006\\
-2.98328875582855	5.00738621111576\\
-2.74548058069763	5.5457768695911\\
-2.48444717307149	6.03839123235924\\
-2.20229858379427	6.48311044821322\\
-1.90123410665799	6.87809280661135\\
-1.58353468812474	7.22155636325303\\
-1.25156052838935	7.51176587237588\\
-0.907745494765833	7.74718733439461\\
-0.554584967057168	7.92662495581357\\
-0.194620615405962	8.04922031186267\\
0.169572435832508	8.1143472468405\\
0.535394186260328	8.12152795950082\\
0.900227241749224	8.07046826719408\\
1.26144841367545	7.96120507451324\\
1.61644712410579	7.7942610021022\\
1.96264806707448	7.57067868796565\\
2.29753142496655	7.29187799317265\\
2.6186441985039	6.95940438672151\\
2.9236015223654	6.57473607222504\\
3.21008419941115	6.1393090213412\\
3.47584259351257	5.6547741889131\\
3.71871299816627	5.1232941004198\\
3.93664172727442	4.54758910213534\\
4.12770297762522	3.93061472196487\\
4.29009991929109	3.27513984996892\\
4.42215590678637	2.58372999227031\\
4.52231729862986	1.85927925584212\\
4.5891774045946	1.10541624138379\\
4.62149671041394	0.325933008355942\\
4.61818817351303	-0.476194372834894\\
4.57829304427835	-1.29824709038394\\
4.5010066391439	-2.13610500986734\\
4.38570763521026	-2.98574878702654\\
4.23190115718114	-3.84558460854962\\
4.03925212865148	-4.71011271778645\\
3.80758344751251	-5.57782034674715\\
3.53685886111262	-6.44396784719875\\
3.22721647058485	-7.30555541868204\\
};
\addlegendentry{sw ($\mathcal{P}$)}

\addplot [color=red, style = thick, dashed]
  table[row sep=crcr]{%
3.22721647058485	-5.17585241577678\\
3.00052204957439	-4.8932417047079\\
2.78629135062245	-4.62233703195288\\
2.58400010897078	-4.36301306257071\\
2.39312968392719	-4.11514512562824\\
2.21316429571315	-3.87873056061739\\
2.04359045287242	-3.65367011866161\\
1.88389557693716	-3.44002221668115\\
1.73356633809088	-3.23776143854693\\
1.59208942351647	-3.04691211068103\\
1.4589489947378	-2.86756167079736\\
1.33362567294576	-2.69977943533446\\
1.21559718910487	-2.54362396637257\\
1.10433787267616	-2.39918727410359\\
0.999317263016812	-2.26658953617311\\
0.899999671617242	-2.14594216442746\\
0.805844935132174	-2.0373319037291\\
0.716309197106727	-1.9408354838816\\
0.630844955445418	-1.85653756600074\\
0.548900716545489	-1.78453025305583\\
0.469921013182036	-1.72489743913777\\
0.393347089819421	-1.6777008510163\\
0.31861793060641	-1.6429788117952\\
0.245171123346099	-1.62075476540886\\
0.172443375955144	-1.61104417789564\\
0.0998709594772396	-1.61385089767333\\
0.0268905379570846	-1.62915362994191\\
-0.047059385767827	-1.65689213051836\\
-0.122536657714023	-1.6969647407115\\
-0.200093920214736	-1.74924237181473\\
-0.280277917861657	-1.81359260524762\\
-0.363629771865317	-1.88989856101493\\
-0.450685660614142	-1.97805805546195\\
-0.541976990742538	-2.07796183419812\\
-0.638029511578578	-2.18946843809367\\
-0.739361811112777	-2.3124018649087\\
-0.846484456722297	-2.44658249480207\\
-0.959900724394516	-2.59186665652374\\
-1.08010827997828	-2.7481493324798\\
-1.20759986586048	-2.91531732886209\\
-1.34286213479827	-3.09321385167168\\
-1.48637486375197	-3.28169076113303\\
-1.63861333725731	-3.48069700795658\\
-1.80005159476271	-3.6902285984051\\
-1.97116120265468	-3.91017998318692\\
-2.15240971784732	-4.14047894273135\\
-2.34426688127544	-4.38129487600855\\
-2.54720271611517	-4.63247110976796\\
-2.76168876798156	-4.89423223611047\\
-2.98820072129521	-5.16653760734458\\
-3.22721647058485	-5.44952638917439\\
};
\addlegendentry{st ($\mathcal{P}$)}

\end{axis}
\end{tikzpicture}%
        \caption*{(b)}
    \end{minipage}
    \caption{Comparison between the solutions of 
$\hat{\mathcal{P}}$ and $\mathcal{P}$ for an average speed of $v_{\text{avg}} = 0.05\sqrt{gl_\circ}$, using a discretization with $51$~points. The period time and the state trajectories are captured accurately. The input provides a qualitatively good approximation. As boundary constraints for the lower-level, we impose~$\hat{\vb}_0$, meaning that the lifted constraints are applied only at the initial time.}
    \label{fig:CGW results}
\end{figure}
The compass-gait walker, shown in Figure~\ref{fig:examples}b, is a bipedal robot consisting of a stance- and a swing leg for which we describe the orientation relative to the vertical with $\theta_{\text{st}}$ and $\theta_{\text{sw}}$ respectively. Furthermore, a motor torque is provided at the hip. The task is to perform a periodic forward motion at a desired average speed $v_{\text{avg}}$, while minimizing control expenditure. Due to the left-right symmetry of the system, we optimize over a single step and obtain a complete stride by flipping the generalized coordinates of stance and swing afterwards. In addition, velocities jump upon the collision of the swing leg with the ground. For the details on the continuous dynamics~\eqref{eq:dynamics} and the discrete jump map~$\Delta:\mathbb{R}^4\to\mathbb{R}^4$, the reader is referred to \cite{Manchester2010}.
The trajectory optimization problem~$\mathcal{P}$ is subject to the MBCs described in Figure~\ref{fig:examples}b.
These MBCs follow a similar structure to the previous examples, incorporating symmetric periodicity and the operating condition at $v_{\text{avg}}$. 
The anchor is chosen at the instance of touch down, when both legs are in contact with the ground. 

Again, the problem was rewritten in the bilevel form~$\hat{\mathcal{P}}$, introduced in Section~\ref{sec:bilevel}.
The boundary conditions were lifted at the initial time with $\hat{\vb}_0$.
To convexify the dynamics, we lifted them by identifying a bilinear Koopman generator surrogate model to apply the presented method. Here, we used a 29-dimensional lift consisting of a mix of trigonometric functions and polynomials of the state. As in the previous example, data was obtained by collecting $45000$ uniformly distributed samples in the state space which were subsequently lifted and used to compute the Lie-derivatives. With this data we solved two least-squares problems to obtain the matrices $\vL_{0}, \, \vL_{1}$ for the bilinear approximation of the Koopman generator.

The solutions of $\hat{\mathcal{P}}$ are compared to those of $\mathcal{P}$ in Figure~\ref{fig:CGW results} for an average speed of $v_{\text{avg}} = 0.05\sqrt{gl_\circ}$.
We observe that we the lifted dynamics are able to fairly accurately compute the optimal state- and input trajectories (with a PCC of $0.9999$ and $0.9987$, respectively) as well as the optimal period $T^\ast$ ($2.2635\, \sqrt{l_\circ/g}$ compared to $2.2518\, \sqrt{l_\circ/g}$).
\section{Discussion}
In this work, we have investigated how the Koopman operator framework can be leveraged to solve non-convex trajectory optimization problems~$\mathcal{P}$ including nonlinear- dynamics and mixed boundary constraints. The problem~$\mathcal{P}$ is reformulated as a bilevel optimization problem~$\hat{\mathcal{P}}$, encompassing both an upper-level ($\Phighhat$) and a lower-level problem ($\Plowhat$). By exploiting this structure, we are able to achieve a convexification of $\Plowhat$. This yields a low-dimensional $\Phighhat$ which can be solved efficiently due to the convexity of $\Plowhat$. 
In addition, an initial guess must only be provided for the boundary values and the terminal time while the initialization of the trajectories is eliminated. The bilevel structure handles the nonlinear boundary constraints well, by including them in the upper-level and not in the lower-level. Furthermore, we have demonstrated the effectiveness of this method using two examples from the domain of periodic trajectory optimization while exploring different formulations of how the boundary constraints are represented in the lower-level.
The results indicate that enforcing one boundary as a hard constraint outperforms soft constraints. For the original problem~$\mathcal{P}$, hard constraints enable more accurate computation of the optimal state, input trajectories, and terminal time.
The presented method is still subject to a number of limiting factors. As demonstrated in~\cite{Otto2024}, linear time-invariant approximations may not be able to fully capture dynamics with products of state and input. While this could potentially pose a problem in cases where control inputs tend to become large, we did not find it to be an issue in our examples where the inputs remained small.
Furthermore, the drift from the lifting manifold presents a challenge, particularly for long-term predictions, as the observables do not span a Koopman-invariant subspace. While re-projecting the lifted state onto the manifold during optimization could address this, it involves the observable functions and thus compromises the problem's convexity.
Future work should explore combining the presented approach with deep learning to jointly learn the Koopman generator and the observables, as proposed by~\cite{Lusch2018} and~\cite{Han2020}. This approach could potentially reduce prediction error while requiring fewer lifting functions.

In our problem definition of $\mathcal{P}$, we made the assumption that the cost must be jointly convex in the state and the input.
This assumption can be potentially relaxed to the requirement that the cost function is jointly convex in the \emph{lifted} state defined by the observables and the input. In addition, we may include nonlinear path constraints in the lower-level as long as it is possible to design the observables such that those constraints define a convex set in the lifted space. This straightforward extension would allow us to solve a larger class of problems.
Finally, we could extend the optimization to include system parameters in the upper-level, similar to the terminal time, forming a co-design problem. To achieve this, it would be essential to treat these parameters as additional states with zero derivatives, ensuring that the Koopman generator does not need to be re-identified in the lower level.

In this work, we have addressed the challenge of incorporating constraints within the Koopman framework by formulating the bilevel optimization problem and investigating three different methods to translate the constraints to the lifted trajectory optimization problem.
\newpage

\acks{This work was funded by the Deutsche Forschungsgemeinschaft (DFG, German Research Foundation) – 501862165. It was further supported through the International Max Planck Research School for Intelligent Systems (IMPRS-IS) for Maximilian Raff.}

\bibliography{literature_bibtex}

\end{document}